\documentclass{amsart}
\usepackage{amsfonts,amssymb,amsthm}
\usepackage[latin9]{inputenc}
\usepackage{mathrsfs}
\usepackage{dcpic,pictex}

\newtheorem{theorem}{Theorem}[section]
\newtheorem{lemma}[theorem]{Lemma}
\newtheorem{corollary}[theorem]{Corollary}
\newtheorem{proposition}[theorem]{Proposition}

\newtheorem{claim}{Claim}

\theoremstyle{definition}
\newtheorem{definition}[theorem]{Definition}
\newtheorem{example}[theorem]{Example}

\theoremstyle{remark}

\numberwithin{equation}{section}

\author[Paola Bonacini]{Paola Bonacini}
\title[Plane section of curves in characteristic $p$]{On the plane section of an integral curve\\ in positive characteristic} 

\address{Dipartimento di Matematica e Informatica, Università degli Studi di Catania, Viale A. Doria 6, 95125 Catania, Italy}
\email{bonacini@dmi.unict.it}

\begin{document}
\begin{abstract}
If $C\subset \mathbb P^3_k$ is an integral curve and $k$ an algebraically closed field of characteristic 0, it is known that the points of the general plane section $C\cap H$ of $C$ are in uniform position. From this it follows easily that the general minimal curve containing $C\cap H$ is irreducible. If char$k=p>0$, the points of $C\cap H$ may not be in uniform position. However, we prove that the general minimal curve containing $C\cap H$ is still irreducible. 
\end{abstract}
\maketitle
\pagestyle{plain}

\section{Introduction}

Let $C$ be an integral curve in $\mathbb P^3_k$, where $k$ is an algebraically closed field of characteristic 0. Harris in \cite{Harris} proves that the points of a general plane section $C\cap H$ of $C$ lie in uniform position, which means that all the subsets of $C\cap H$ of the same cardinality have the same Hilbert function. A set of points with the Uniform Position Property has some interesting properties, as the fact that the general minimal curve containing it is irreducible.
 
In positive characteristic the Uniform Position Property is not satisfied for all the generic plane section of integral curves. Indeed in \cite[Example 1.2]{Rath} Rathmann gives an example of a complete intersection integral curve, whose generic plane section is not in uniform position. However the general minimal curve containing it is still irreducible. 

In \S 3 we prove that, independently of the characteristic of the base field, the general minimal curve containing the generic plane section of an integral curve is always irreducible. From it also follows easily that the generic plane section of an integral curve has the Hilbert function of decreasing type, fact well known in characteristic 0 (see \cite{Harris} for example) and proved in positive characteristic in \cite{B}.

\section{Preliminaries and notation}

Let us denote $X$ a finite set of points in $\mathbb P^2_k$, where $k$ is an algebraically closed field. Let us recall the following:
\begin{definition}
We say that $X$ has the \emph{Uniform Position Property} (briefly UPP) if, for any $n\le \deg X$, all the subsets of $n$ points have the same Hilbert function.
\end{definition}
If this happens, then (see \cite{GMP}) the Hilbert function of a set $Z\subset X$ of $n$ points is the following truncated function:
\begin{equation}  \label{E:1}
H(Z,i)=\min\{H(X,i),n\}\ \forall i.
\end{equation}

One of the main consequences of the Uniform Position Property is the following:
\begin{proposition}[{\cite[p.198]{Harris}}]  \label{P:1}
If $X$ is a set of points in $\mathbb P^2$ in uniform position and $X'\subset X$ is such that $h^1(\mathscr I_{X'}(l))\ne 0$, then $h^0(\mathscr I_X(l))=h^0(\mathscr I_{X'}(l))$. So every curve of degree $l$ containing $X'$ contains $X$. In particular this holds when $\deg X'=\binom{l+2}{2}$.
\begin{proof} 
If $h^1(\mathscr I_{X'}(l))\ne 0$, then:
\[H(X',l)=h^0\left(\mathscr O_{\mathbb P^2}(l)\right)-h^0\left(\mathscr I_{X'}(l)\right)< \deg X'\]
and so, by \eqref{E:1}, $H(X',l)=H(X,l)$, which implies that:
\[h^0\left(\mathscr I_{X'}(l)\right)=h^0\left(\mathscr I_X(l)\right).\]

If $\deg X'=\binom{l+2}{2}$ and $h^0(\mathscr I_{X'}(l))\ne 0$, one has: 
\[\chi(\mathscr I_X(l))=\chi(\mathscr O_{\mathbb P^2}(l))-\chi(\mathscr O_X(l))=0,\]
so that:
\[h^1(\mathscr I_{X'}(l))=h^0(\mathscr I_{X'}(l))\ne 0.\] 
If $h^0(\mathscr I_{X'}(l))=0$, then it must be $h^0(\mathscr I_X(l))=0$.
\end{proof}
\end{proposition}

An important corollary of this fact is:
\begin{proposition}[{\cite[Remark 1.2(ii)]{MR1}},{\cite[Theorem 3.4]{GM}}] \label{P:2}
Let $X\subset \mathbb P^2$ be a reduced scheme in uniform position such that:
\[\deg X=n=\binom{d+2}{2}+h\]
with $0\le h\le d+1$. If:
\[g=\inf\left\{i\in \mathbb N|H^0\left(\mathscr I_X(i)\right)\ne 0\right\}\]
then:
\begin{enumerate}
	\item if either $g\le d$ or $g=d+1$ and $h\ge 2$, every element of $H^0(\mathscr I_X(g))$ is irreducible;
	\item if $g=d+1$ and $h=0,1$, the generic element of $H^0(\mathscr I_X(g))$ is irreducible. In particular there exists a minimal smooth curve containing $X$.
\end{enumerate}
\begin{proof}
We follow the proof in \cite{MR1}. First of all note that $h^0(\mathscr I_X(d+1))\ne 0$, so that $g\le d+1$. Indeed:
\[\binom{d+1+2}{2}-\deg X=\binom{d+3}{2}-\binom{d+2}{2}-h=d+2-h\ge 1.\]

If $F=F_1F_2\in H^0(\mathscr I_X(g))$, with $\deg F_1=r_1>0$, $\deg F_2=r_2>0$ and $r_1+r_2=g$, then, by Proposition \ref{P:1}, $F_1$ (resp. $F_2$) can contain at most $\binom{r_1+2}{2}-1$ (resp.  $\binom{r_2+2}{2}-1$) points of $X$. So:
\[n\le \binom{r_1+2}{2}-1+\binom{r_2+2}{2}-1\]
from which it follows that:
\[d^2+3d+2h\le {r_1}^2+{r_2}^2+3g-2.\]
But ${r_1}^2+{r_2}^2=g^2-2r_1r_2$, where:
\[r_1r_2=r_1(g-r_1)\ge g-1\]
and so:
\[d^2+3d+2h\le g^2+g.\]
If $g\le d$, then we get a contradiction for any $h$. If $g=d+1$, then we conclude that $h\le 1$. In this case, not every curve containing $X$ of degree $g=d+1$ is irreducible, but using \cite[Theorem 3.5]{MR} we see that this is true for the generic one, because there exists one which is nonsingular. 
\end{proof} 
\end{proposition}

The most important result relating UPP and plane section of curves is the following:
\begin{theorem}[Uniform Position Lemma, {\cite[\S 2]{Harris}}] \label{T:1}
If $C\subset \mathbb P^3_k$ is an integral curve and \emph{char}$k=0$, then the points of the generic plane section of $C$ are in uniform position.
\end{theorem}

\begin{corollary}
If $C\subset \mathbb P^3_k$ is an integral curve and \emph{char}$k=0$, the general minimal curve containing the generic plane section of $C$ is irreducible.
\begin{proof}
Just apply Proposition \ref{P:2} and Theorem \ref{T:1}.
\end{proof}
\end{corollary}

\begin{example}[see {\cite[Example 1.2]{Rath}}]
Let $C$ be the curve defined as the complete intersection of ${x_0}^q-x_1{x_3}^{q-1}={x_1}^q-x_2{x_3}^{q-1}=0$, where $q=p^f$, for some $f>0$. On the open affine set defined by $x_3\ne 0$ the curve $C$ is described by:
\[x_0=t,\ x_1=t^q,\ x_2=t^{q^2}\]
so that $C$ is rational and integral. However the general plane section of $C$ looks like a 2-dimensional vector space over a field with $q$ elements. Indeed, if $P_0,P_1,P_2\in C$ are linearly independent and $H$ is the plane containing the three points, then the points of $C\cap H$ are given by:
\[P=P_0+\lambda_1(P_1-P_0)+\lambda_2(P_2-P_0)\]
where $\lambda_i^q=\lambda_i$ for $i=1,2$. Note also that given two points $P_0,P_1\in C$ the line $P_0P_1$ meets $C$ in the points given by:
\[P=P_0+\mu(P_1-P_0)\]
where $\mu^q=\mu$. So in every generic plane section of $C$ there are, for $q\ge 4$, at least 3 collinear points and 3 linearly independent. It means that the generic plane section of $C$ has not the UPP. 

Since $C$ is a complete intersection of two surfaces of degree $q$, $C\cap H$ is a complete intersection too and the minimal curves containing $C\cap H$ have degree $q$. Since $C$ is integral, these surfaces are integral and the generic plane section of each of them is an integral curve in $H$. So, in this case, every minimal curve containing the generic plane section of $C$ is irreducible, because every minimal surface containing $C$ is irreducible.
\end{example}

\section{The fibres of morphisms}

In this section we recall some results that will be useful in the sequel. 

\begin{lemma}  \label{L:2}
Let $Y$ be a locally noetherian and integral scheme and $f:X\rightarrow Y$ a morphism of finite type. Then there exists an open nonempty subset $U$ of $Y$ such that $f|_{f^{-1}(U)}:f^{-1}(U)\rightarrow U$ is flat.
\begin{proof}
This follows from \cite[Th\'eor\`eme 6.9.1]{EGA}, with $\mathscr F=\mathscr O_X$.
\end{proof}
\end{lemma}

\begin{lemma}[{\cite[Corollaire 3.3.5]{EGA}}]  \label{L:1}
Let $X$ and $Y$ be locally noetherian schemes and $f:X\rightarrow Y$ be a flat morphism. If $Y$ is reduced in the points of $f(X)$ and if $f^{-1}(y)$ is a $k(y)$-scheme reduced for every $y\in f(X)$, then $X$ is reduced.
\end{lemma}

Another result that will be useful in the proof of the main result of this paper is the following:
\begin{theorem}   \label{T:4}
Let $f:X\rightarrow Y$ be a dominant morphism from an irreducible algebraic variety $X$ to an algebraic variety $Y$. Let $U\subseteq Y$ and $V\subseteq f^{-1}(U)\subset X$ be affine open subsets and $u_1,\dots,u_r\in k[V]$ algebraically independent over $k(U)=k(Y)$ such that for some variable $T$:
\[ k[V]\cong k[U][u_1,\dots, u_r][T]/(F)\]
with $F\in k[U][u_1,\dots,u_r][T]$ irreducible polynomial. If $k(Y)$ is algebraically closed in $k(X)$, there exists a nonempty open subset $W\subset V$ such that for every $y\in W$ the fibre $f^{-1}(y)$ is geometrically integral.
\begin{proof}
The proof works as in \cite[Theorem 1, p.139]{Sha}.
\end{proof}
\end{theorem}

\section{Curves of least degree}

In this section we consider an integral curve $C\subset\mathbb P^3_k$, where $k$ is an algebraically closed field of any characteristic. Let $\{\underline x_i\}$ be a system of coordinates over $\mathbb P^3$ and $\{\underline t_j\}$ a system of coordinates over the dual space ${\mathbb P^3}^{\vee}$. Denote by $\mathscr I_C$ the ideal sheaf of $C$ in $\mathbb P^3$ and let us consider the incidence variety $M\subset {\mathbb P^3}^{\vee}\times \mathbb P^3$ determined by:
\[h=\sum_{i=0}^3t_ix_i=0.\]
Let us consider the two projections:
\[p:M\rightarrow \mathbb P^3,\]
\[g:M\rightarrow {\mathbb P^3}^{\vee}\]
and the following module over $\mathscr O_M$:
\[\mathscr I(m,n)=g^{\star}\left(\mathscr O_{{\mathbb P^3}^{\vee}}(m)\right)\otimes_{\mathscr O_M} p^{\star}\left(\mathscr I_C(n)\right).\]
If $T={p}^{-1}(C)$, then it is easy to see that $T$ is integral. Indeed, every fibre of the projection $p|_T:T\rightarrow C$ is integral and equidimensional of dimension 2, being isomorphic to a projective plane. So $T$ is irreducible and $p|_T$ is flat, by \cite[Ch. III, Theorem 9.9]{Hart} (note that $T\subset {\mathbb P^3}^{\vee}\times C$). By Lemma \ref{L:1} $T$ is reduced. It is also possible to show that $\mathscr I=\mathscr I_T$, the ideal sheaf of $T$ in $M$ (see \cite{B}).

\begin{theorem}   \label{T:3}
If $X\subset \mathbb P^2$ is a generic plane section of an integral curve $C\subset \mathbb P^3$, then the generic curve in $\mathbb P^2$ of minimal degree containing $X$ is irreducible. In particular, if there exists only one curve of minimal degree containing $X$, such a curve is irreducible.
\begin{proof}
Let us denote by $s$ the minimal degree of a plane curve containing $X$. Then there exists $\alpha\gg 0$ such that $H^0(\mathscr I(\alpha,s))\ne 0$ (see \cite{GP2} and \cite[Proposition 3.2]{B}). So we can take a general form $F\in H^0(M,\mathscr I(\alpha,s))$, with $\alpha$ minimal, in such a way that it determines a hypersurface $S$ in the incidence variety $M$. Since $\mathscr I=\mathscr I_T$, $T\subset S$. $T$ is integral and so $S$ is integral too, by the minimality of $s$ and $\alpha$.

Since the projection $T\rightarrow {\mathbb P^3}^{\vee}$ is dominant, $S\rightarrow {\mathbb P^3}^{\vee}$ is dominant too and, as in \cite[Proposition 3.2]{B}, we see that for a generic $P_H\in {\mathbb P^3}^{\vee}$ corresponding to a generic plane in $\mathbb P^3$, $p\,(S\cap g^{-1}(P_H))$ is a curve of degree $s$ containing $C\cap H$. 

Note that by the K{\"u}nneth formula (see \cite[Ch. VI, Corollary 8.13]{Milne}) we get:
\[H^1\left(\mathscr O_{{\mathbb P^3}^{\vee}\times \mathbb P^3}(m,n)\right)=0\ \forall\, m,n.\]
This implies that we have the following surjective map:
\[H^0\left(\mathscr O_{{\mathbb P^3}^{\vee}\times \mathbb P^3}(m,n)\right)\rightarrow H^0\left(\mathscr O_M(m,n)\right)\rightarrow 0\ \forall\, m,n.\] 
So, since $\mathscr I\subset \mathscr O_M$, the form $F\in H^0(\mathscr I(\alpha,s))$ defining $S$ can be lifted to a homogeneous polynomial $f\in k[\underline t;\underline x]$ of bidegree $(\alpha,s)$. So $S$ in ${\mathbb P^3}^{\vee}\times \mathbb P^3$ is the complete intersection determined by:
\begin{equation*}  
\begin{split}
f(\underline t;\underline x)=0\\
\sum_{i=0}^3t_ix_i=0.
\end{split}
\end{equation*}
Note that restricting to the open subsets $U\subset {\mathbb P^3}^{\vee}$, given by $t_3\ne 0$, and $V\subset S$, given by $t_3\ne 0, x_2\ne 0$, we get:
\[k[U]=k[t_0,t_1,t_2]\]
\begin{equation} \label{S:18}
k[V]=k[t_0,t_1,t_2,x_0,x_1]/(F)
\end{equation} 
where $F(t_0,t_1,t_2,x_0,x_1)=f(t_0,t_1,t_2,1;x_0,x_1,1,t_0x_0+t_1x_1+t_2)$.

We want to show that the generic curve of degree $s$ containing $X$ is irreducible. So we show that the generic fibre of the projection $g_S:S\rightarrow {\mathbb P^3}^{\vee}$ is irreducible. By Theorem \ref{T:4} we see that it's sufficient to show that $k({\mathbb P^3}^{\vee})$ is algebraically closed in $k(S)$, since by \eqref{S:18} we see that in this case all the conditions are satisfied. 

Let us consider the algebraic closure $K$ of ${g_S}^{\star}k({\mathbb P^3}^{\vee})$ in $k(S)$. We want to show that $K={g_S}^{\star}k({\mathbb P^3}^{\vee})$. Of course $K$ is finitely generated over $k$ and so it corresponds to an irreducible variety $Z$ such that $k(Z)=K$. Since $K$ is algebraic over ${g_S}^{\star}k({\mathbb P^3}^{\vee})$, we see that $\dim Z=3$ and that we can factor $g_S$ as a composition of rational maps:
\begin{equation} \label{S:3}
S\stackrel{\pi_1}{\dashrightarrow} Z\stackrel{\pi_2}{\dashrightarrow}{\mathbb P^3}^{\vee}
\end{equation}  
where $\pi_1$ is dominant and has irreducible fibres (by Theorem \ref{T:4}) and $\pi_2$ is generically finite.

Let us take an open subset $U\subset {\mathbb P^3}^{\vee}$, whose points correspond to planes that are general for $C$, and $V={\pi_2}^{-1}(U)\subset Z$ and $W={g_S}^{-1}(U)\subset S$, so that the following:
\[W\stackrel{\pi_1}{\rightarrow} V\stackrel{\pi_2}{\rightarrow} U\]
is a factorisation of $g_U:W\rightarrow U$ and the following conditions are satisfied: 
\[g_U\ \text{is flat,}\]
\[\text{the projection\ } T\cap(g^{-1}(U))\rightarrow U\ \text{is flat, with}\ g:M\rightarrow {\mathbb P^3}^{\vee},\]
\begin{equation}  \label{S:5}
{g_U}^{-1}(P_H)\text{ is reduced and } \dim {g_U}^{-1}(P_H)=1\ \text{for any}\ P_H\in U
\end{equation}
and:
\begin{equation} \label{S:14}
{\pi_2}^{-1}(L\cap U)\ \text{is irreducible}
\end{equation} 
for any $L\subset {\mathbb P^3}^{\vee}$ plane such that $L\cap U\ne \emptyset$. Note that $U\ne \emptyset$ by Lemma \ref{L:2} and by \cite[I, Theorem 6.10(3)]{Jou}. 

We will show that $\pi_2:V\rightarrow U$ is an isomorphism.

Let us consider the fibred product $W\times_U V$ and the projection:
\[\phi:W\times_U V\rightarrow V.\]
\begin{claim}
$W\times_UV$ is equidimensional of dimension 4 and is reduced.
\begin{proof}[proof of Claim 1]
By \eqref{S:5} any point $Q_0\in V$ is such that $\pi_2(Q_0)=P_H\in U\subset {\mathbb P^3}^{\vee}$ corresponds to a plane $H$ and so:
\begin{equation}  \label{S:10}
\phi^{-1}(Q_0)=(W\times_UV)\times_V \operatorname{Spec}k(Q_0)\cong W\times_U\operatorname{Spec}k(P_H)={g_U}^{-1}(P_H).
\end{equation}
So $\phi^{-1}(Q_0)$ is isomorphic to a curve in $H$ of degree $s$ containing $C\cap H$. So all the fibres of $\phi$ are equidimensional of dimension 1 and so we see that $W\times_UV$ is equidimensional of dimension 4. $\phi^{-1}(Q_0)$ is also reduced for any $Q_0\in V$ and $\phi$ is flat, by base change, being $g_U$ flat by construction. So using Lemma \ref{L:1} we see that $W\times_UV$ is reduced.
\renewcommand{\qedsymbol}{}
\end{proof}
\end{claim}

\
Let us now consider $g:M\rightarrow {\mathbb P^3}^{\vee}$ and $T=p^{-1}(C)$. Since $g^{-1}(U)$ is an open subset of $M$ and $T\subset M$, then $T\cap(g^{-1}(U))$ is dense in $T$. So we can consider $T'=(T\cap(g^{-1}(U)))\times_UV$. 

\begin{claim}
$T'$ is integral.
\begin{proof}[proof of Claim 2]
Let us consider the projection:
\[\psi:T'\rightarrow C\]
so that $\psi(T')$ is an open subset of $C$. Taken any $P_0\in \psi(T')$, we see that:
\begin{equation}  \label{S:9}
\psi^{-1}(P_0)\cong (L_{P_0}\cap U)\times_UV\cong {\pi_2}^{-1}(L_{P_0}\cap U)
\end{equation}
where $L_{P_0}\subset \mathbb P^3$ is the linear subvariety of the planes containing $P_0$ and has dimension 2. By the definition of $T'$ we see that $L_{P_0}\cap U\ne \emptyset$ and by \eqref{S:14} ${\pi_2}^{-1}(L_{P_0}\cap U)$ is irreducible for any $P_0\in \psi(T')$. So from \eqref{S:9} we conclude that $T'$ is irreducible. 

Let us consider now $\phi_{T'}:T'\rightarrow V$. This is a surjective morphism such that:
\begin{equation}  \label{S:12}
{\phi_{T'}}^{-1}(Q_0)=T'\times_V\operatorname{Spec}k(Q_0)\cong T\times_U\operatorname{Spec}k(P_H)\cong C\cap H.
\end{equation}
Since by construction $T\cap(g^{-1}(U))\rightarrow U$ is flat, by base change $\phi_{T'}$ is flat too and so by Lemma \ref{L:1} we see that $T'$ is reduced. So the claim is proved.
\renewcommand{\qedsymbol}{}
\end{proof}
\end{claim}

\begin{claim}
$W\times_UV$ is integral. 
\begin{proof}[proof of Claim 3]
Since $\dim W\times_UV=4$ and $W\times_UV$ is reduced, there exists an integral component $\tilde W$ of $W\times_UV$, with $\dim \tilde W=4$, such that $T'\subset \tilde W$. Since $\phi_{T'}$ is surjective, then: 
\[\varphi:\tilde W\rightarrow V\] 
is surjective too. So, taken any $Q_0\in V$, if $\pi_2(Q_0)=P_H$, we have:
\begin{equation}  \label{S:15}
\dim  {\varphi}^{-1}(Q_0)\ge 1
\end{equation}
and
\begin{equation}   \label{S:13}
C\cap H\cong{\phi_{T'}}^{-1}(Q_0)\subset {\varphi}^{-1}(Q_0)\subset {\phi}^{-1}(Q_0)={g_U}^{-1}(P_H)
\end{equation}
by \eqref{S:10} and \eqref{S:12}. By \eqref{S:5} ${\phi}^{-1}(Q_0)={g_U}^{-1}(P_H)$ is equidimensional of dimension 1 for any $Q_0\in V$ and we see by \eqref{S:15} that ${\varphi}^{-1}(Q_0)\subset {g_U}^{-1}(P_H)$ is equidimensional of dimension 1 for any $Q_0\in V$. So ${\varphi}^{-1}(Q_0)$ is a component of the curve ${g_U}^{-1}(P_H)$ of degree $s$ containing $C\cap H$ by \eqref{S:13}. By the minimality of $s$ $\varphi^{-1}(Q_0)\cong {g_U}^{-1}(P_H)={\phi}^{-1}(Q_0)$ for any $Q_0\in V$. So $\tilde W=W\times_UV$, which means that $W\times_UV$ is reduced and irreducible.
\renewcommand{\qedsymbol}{}
\end{proof}
\end{claim}

\begin{claim}
$W\times_UV=W$.
\begin{proof}[proof of Claim 4]
Let us now consider the projection $\psi:W\times_UV\rightarrow W$ and take any $P_0\in W$. Then:
\[{\psi}^{-1}(P_0)=\operatorname{Spec}k(P_0)\times_UV\cong {\pi_2}^{-1}(g(P_0)).\]
So $\psi$ is a finite morphism and: 
\begin{equation}  \label{S:11}
\psi^{\star}k(W)\subset k(W\times_UV)
\end{equation}
is an algebraic field extension. Obviously we have the following commutative diagram:
\[
\begin{vcenter}
{\begindc{\commdiag}
\obj(1,1)[W]{$W$}
\obj(3,1)[U]{$U$}
\obj(1,3)[W1]{$W$}
\obj(3,3)[V]{$V$}
\mor{W}{U}{$g_U$}[\atright,\solidarrow]
\mor{V}{U}{$\pi_2$}
\mor{W1}{W}{$i$}[\atright,\solidarrow]
\mor{W1}{V}{$\pi_1$}
\enddc}
\end{vcenter}
\]
and so, by the properties of the fibred product, there exists a morphism $\tau:W\rightarrow W\times_UV$ that factors the identity on $W$:
\begin{equation}  \label{S:6}
\begin{vcenter}
{\begindc{\commdiag}
\obj(1,3)[W]{$W$}
\obj(5,3)[W1]{$W\times_UV$}
\obj(3,1)[W2]{$W$}
\mor{W}{W1}{$\tau$}[\atleft,\dashArrow]
\mor{W1}{W2}{$\psi$}
\mor{W}{W2}{$i$}[\atright,\solidarrow]
\enddc}
\end{vcenter}
\end{equation}
So, being $W\times_UV$ integral, $\tau$ must be surjective and so we get the inclusion ${\tau}^{\star}k(W\times_UV)\subset k(W)$. Applying ${\tau}^{\star}$ to \eqref{S:11} we get:
\[k(W)={\tau}^{\star}{\psi}^{\star}k(W)\subset {\tau}^{\star}k(W\times_UV)\subset k(W)\]
so that ${\tau}^{\star}k(W\times_UV)=k(W)$. It means that $\tau$ is an isomorphism and from \eqref{S:6} we see that $\psi$ is an isomorphism.
\renewcommand{\qedsymbol}{}
\end{proof}
\end{claim}

\
From $W=W\times_UV$ it follows immediately that $U=V$ and in this way the theorem is proved.
\end{proof}
\end{theorem}

Now we show that from this result it follows easily that the Hilbert function of the generic plane section of $C$ is of decreasing type. Let $X=C\cap H$ be a general plane section of $C$ and $H(X,\_ \ )$ the Hilbert function of $X$ in $H$. The first difference of the Hilbert function is:
\[\Delta H(X,i)=H(X,i)-H(X,i-1)\ \forall i.\]
It is known \cite{GMP} that there exists $a_1\le a_2\le t$ such that:
\[
\Delta H(X,i)=\begin{cases}
    i+1                   & \text{for $i=0,\dots,a_1-1$}\\
    a_1                   & \text{for $i=a_1,\dots,a_2-1$}\\
    <a_1                  & \text{for $i=a_2$}\\
    \text{non increasing} & \text{for $i=a_2,\dots,t$}\\
    0                     & \text{for $i>t$.}\\
\end{cases}    
\]    
\begin{definition}
We say that $X$ has \emph{the Hilbert function of decreasing type} if for $a_2\le i<j<t$ we have $\Delta H(X,i)>\Delta H(X,j)$.
\end{definition}
It is well known that the Hilbert function of the generic plane section of an integral curve is of decreasing type (see \cite{Harris}, \cite{MR2}, \cite[Lemme 3.2]{GP} for char$k=0$ and \cite[Theorem 3.3]{B} for any characteristic). This fact is direct consequence of Theorem \ref{T:3}:
\begin{corollary}
Let $C\subset \mathbb P^3$ be an integral curve and let $X=C\cap H$ be a general plane section of $C$. Then the Hilbert function of $X$ is of decreasing type.
\begin{proof}
Let $c_i=\Delta H(X,i)$ for each $i$. Let us suppose that for some $s\ge a_1$ we have $c_{s-1}=c_s$. From \cite[Theorem 2.9]{MR1} we see that the curves in $H$ of degree $s$ containing $X$ have a common factor $D$ of degree $c_s$. In particular $D$ is a common factor of the minimal curves containing $X$. However from Theorem \ref{T:3} the general minimal curve is irreducible. So $D$ is irreducible and is the only minimal curve containing $X$. Since $\deg D=c_s$, we have $c_s=a_1$ and so the Hilbert function of $X$ is of decreasing type.
\end{proof}
\end{corollary}  

\subsection*{Acknowledgements}
I would like to thank my advisor Rosario Strano for all his useful suggestions and Riccardo Re for the many helpful conversations.

\bibliographystyle{amsplain}

\end{document}